\documentclass[]{article}
\usepackage[dvipsnames]{xcolor}

\usepackage[numbers,sort&compress]{natbib}

\usepackage[utf8]{inputenc}
\usepackage{lmodern}
\usepackage{mathtools} \usepackage{amssymb}
\usepackage{amsthm}

\newtheorem{theorem}{Theorem}[section]
\newtheorem{proposition}[theorem]{Proposition}
\newtheorem{corollary}[theorem]{Corollary}
\theoremstyle{definition}
\newtheorem{definition}[theorem]{Definition}
\newtheorem{remark}[theorem]{Remark}
\newtheorem{example}{Example}[section]

\DeclareMathOperator{\cl}{cl}
\DeclareMathOperator{\bd}{bd}
\DeclareMathOperator{\conv}{conv}
\DeclareMathOperator{\cone}{cone}

\DeclareMathOperator{\interior}{int}
\DeclareMathOperator{\vertices}{vert}

\DeclareMathOperator{\recc}{0^{\infty}}
\DeclareMathOperator{\tr}{trace}

\newcommand{\polar}[1]{{#1}^{\circ}} 
\newcommand{\R}{\mathbb{R}}
\renewcommand{\S}{\mathcal{S}}
\newcommand{\B}{\mathcal{B}}
\newcommand{\T}{^{\mathsf{T}}}
\AtBeginDocument{\renewcommand{\O}{\mathcal{O}}}
\newcommand{\I}{\mathcal{I}}
\newcommand{\set}[1]{\left\lbrace {#1} \right\rbrace}
\newcommand{\norm}[1]{\left\lVert {#1} \right\rVert}
\newcommand{\haus}[2]{d_{\mathsf{H}}\left({#1},{#2}\right)} 
\newcommand{\thaus}[2]{\bar{d}_{\mathsf{H}}(#1,#2)} 
\newcommand{\mgeq}{\succcurlyeq} 
\newcommand{\mg}{\succ} 
\newcommand{\A}{\mathcal{A}}

\newcommand{\decoRule}{\rule{.8\textwidth}{.4pt}}

\usepackage{xifthen}
\AtBeginDocument{\renewcommand{\vector}[2][4]{\ifthenelse{#1=2}{\begin{pmatrix}{#2}_1\\{#2}_2\end{pmatrix}}{
    \ifthenelse{#1=3}{\begin{pmatrix}{#2}_1\\{#2}_2\\{#2}_3\end{pmatrix}}{
      \begin{pmatrix}{#2}_1\\{#2}_2\\\vdots\\{#2}_n\end{pmatrix}}}}}
\AtBeginDocument{\renewcommand{\matrix}[1]{\begin{pmatrix} #1 \end{pmatrix}}}
\DeclareRobustCommand\varC{C}
\DeclareRobustCommand\varw{w}
\DeclareRobustCommand\varv{c}
\DeclareRobustCommand\vard{d}
\newenvironment{optimizationproblem}[2]{\begin{equation}\tag{#1}\label{#2}\begin{aligned}}{\end{aligned}\end{equation}}

\usepackage{hyperref}
\usepackage{color}
\usepackage{enumitem}
\usepackage{subcaption}
\usepackage{booktabs}
\usepackage{makecell}
\usepackage{diagbox}
\usepackage[misc]{ifsym}
\usepackage[linesnumbered,algoruled,vlined]{algorithm2e}
\usepackage{ifthen}
\newcommand{\isassigned}{\leftarrow}

\newcommand{\figpath}{figures/}

\begin{document}

\title{A polyhedral approximation algorithm for recession cones of
  spectrahedral shadows}

\author{Daniel Dörfler \\
  Friedrich Schiller University Jena, Germany
  \and Andreas Löhne \\
  Friedrich Schiller University Jena, Germany}

\date{\today}

\maketitle

\begin{abstract}
  The intersection of an affine subspace with the cone of positive
  semidefinite matrices is called a spectrahedron. An orthogonal
  projection thereof is called a spectrahedral shadow or projected
  spectrahedron. Spectrahedra and their projections can be seen as a
  generalization of polyhedra. This article is concerned with the
  problem of approximating the recession cones of spectrahedra and
  spectrahedral shadows via polyhedral cones. We present two iterative
  algorithms to compute outer and inner approximations to within an
  arbitrary prescribed accuracy. The first algorithm is tailored to
  spectrahedra and is derived from polyhedral approximation algorithms
  for compact convex sets and relies on the fact, that an algebraic
  description of the recession cone is available. The second algorithm
  is designed for projected spectrahedra and does not require an
  algebraic description of the recession cone, which is in general more
  difficult to obtain. We prove correctness and finiteness of both
  algorithms and provide numerical examples.
\end{abstract}

\section{Introduction}
Polyhedral approximation is a common technique used in mathematical
programming. It is, e.g., used in algorithms that approximate the
feasible region of a convex optimization problem by a sequence of
polyhedra in order to obtain an approximate solution. Early
publications of these algorithms include, amongs others, Cheney's and
Goldstein's Newton method for convex programming \cite{CG59}, Kelley's
cutting-plane method \cite{Kel60}, and Veinott's supporting hyperplane
method \cite{Vei67}. These algorithms have since been improved and
extended in various directions and different fields, such as global
optimization \cite{Tuy83, TTB83}, geometry \cite{LL92}, solution
concepts for multiple objective optimization problems \cite{RW05,
HL10, ESS11, LRU14}, and algorithms for mixed-integer convex
optimization problems \cite{DG86, WP95, KLW16}. Kamenev unites the
approximation ideas in the above mentioned works into a family of
iterative algorithms called augmenting and cutting schemes
\cite{Kam92} and analyzes their efficiency and convergence properties
in \cite{Kam93, Kam94}.

Interest in polyhedral approximations arise from the fact that, while
polyhedra do not have to be finite sets, they can be represented in a
finite manner and many set calculus operations can be performed with
them straightforwardly, see \cite{CLW19}. This makes them particularly
suitable for computer-aided treatment and easier to work with than
with more general convex sets.

The algorithms mentioned so far are tailored to compact sets, because
unbounded sets can only be approximated by polyhedra in the Hausdorff
distance if they satisfy restrictive properties as investigated in
\cite{NR95}. However, unbounded sets arise naturally in
applications. In the theory of convex vector optimization, for
example, there exists a solution concept, that is based on finding a
polyhedral approximation of the so-called extended feasible image of
the problem. This set is unbounded by its construction, but it is
advantageous to work with it rather than with the possibly bounded
feasible image, because its boundary satisfies certain minimality
properties, see \cite{Loe11, ESS11}. The recession cone of a set
describes its asymptotic behaviour and is therefore a crucial
characteristic in the study of unbounded sets. In order to properly
deal with unboundedness in the framework of polyhedral approximation,
the recession cone of the given set should be approximated in some
sense. There are several ways to define metrics on the space of closed
and convex cones, which are summarized in the survey article
\cite{IS10}.

Spectrahedra are intersections of affine subspaces with the cone of
positive semidefinite matrices and can be defined by affine matrix
inequalities. Projections of spectrahedra are called projected
spectrahedra, spectrahedral shadows or semidefinitely-representable
sets and arise as the feasible regions of semidefinite
programs. Hence, they are ubiquitous in modern optimization. They can
be viewed as a generalization of polyhedra, as they contain the class
of polyhedral sets as a proper subclass and as they are closed under
many operations under which polyhedra are also closed, see
e.g. \cite{Net11}. One aspect in which spectrahedra and their
projections differ from each other is the difficulty in expressing
their recession cones. While the recession cone of a spectrahedron can
be described by the linear part of the defining matrix inequality,
there is no straightforward way to describe the recession cone of a
projected spectrahedron from its data.

This manuscript is concerned with the task of computing polyhedral
approximations of the recession cones of specrahedra and spectrahedral
shadows. To the best of the authors knowledge, there do not exist any
algorithms in the literature tailored to this problem.  We present two
iterative approximations algorithms for recession cones. The first one
is designed for spectrahedra and is derived from augmenting and cutting
schemes. A key requirement for this algorithm is that an algebraic
description of the recession cone is known. The second algorithm does
not require prior knowledge about the whole recession cone, but only a
single element must be known. It is suitable for the approximation of
the recession cone of a spectrahedral shadow. Both algorithms compute
polyhedral inner and outer approximations to within a prescribed
accuracy. We show that the algorithms are correct and finite, given
that the input set contains no lines and has a solid recession cone.
The algorithms are tested and compared on three examples.

The article is structured as follows. In Section \ref{section2} we
introduce the necessary notation and definitions. Section
\ref{section3} is devoted to the polyhedral approximation of recession
cones of spectrahedra. After reciting some results about convex cones
and presenting first results about relevant semidefinite optimization
problems, the first algorithm is presented. Thereafter, its
correctness and finiteness is proved. The algorithm for the
approximation of recession cones of spectrahedral shadows is provided
in Section \ref{section4}. Again, we prove correctness and
finiteness. In the final section, we test both algorithm on examples.

\section{Preliminaries}\label{section2}
Given a set $C \subseteq \R^n$, we denote by $\conv C$, $\cone C$,
$\interior C$, $\bd C$ and $\cl C$ the \textit{convex hull}, \textit{conical
hull}, \textit{interior}, \textit{boundary}, and \textit{closure} of $C$,
respectively. The \textit{Minkowski addition} between two sets $C_1$
and $C_2 \subseteq \R^n$ is defined as $C_1+C_2 := \set{c_1+c_2 \mid
c_1 \in C_1, c_2 \in C_2}.$ A set $K$ is called a \textit{cone}, if $x
\in K$ implies $\alpha x \in K$ for every $\alpha \geq 0$. A cone $K$
is called \textit{pointed}, if $K \cap (-K) = \set{0}$. The set
\[
\set{d \in \R^n \mid x \in C \Rightarrow x+d \in C},
\]
denoted by $\recc C$, is called the \textit{recession cone} of the
convex set $C$. A closed convex set is called \textit{line-free}, if its
recession cone is pointed. The \textit{polar cone} of a convex cone $K
\subseteq \R^n$ is defined as $\set{x \in \R^n \mid y \in K
\Rightarrow x\T y \leq 0}$ and denoted by $K^{\circ}$. A
\textit{polyhedron} $P \subseteq \R^n$ is defined as the intersection
of finitely many closed halfspaces, i.e. $P = \set{x \in \R^n \mid Ax
\leq b}$ for some $A \in \R^{m \times n}$, $b \in \R^m$, and $m \in
\mathbb{N}$. The tuple $(A,b)$ is called a \textit{H-representation}
or halfspace representation of $P$. Alternatively, $P$ can also be
expressed as the Minkowski sum of the convex hull of finitely many
points and the conical hull of finitely many directions, i.e. there
exist $k, \ell \in \mathbb{N}$, points $v_1, \dots, v_k \in \R^n$, and
unit directions $d_1, \dots, d_{\ell} \in \R^n$, such that $P = \conv
\set{v_i \mid i=1,\dots,k} + \cone \set{d_j \mid j=1,\dots,\ell}$. In
this connection, we use the convention $\conv \emptyset = \emptyset$
and $\cone \emptyset = \set{0}$. Given the matrices $V \in \R^{n
\times k}$ and $D \in \R^{n \times \ell}$, where the $i$-th column of
$V$ is $v_i$ and the $j$-th column of $D$ is $d_j$, we call the tuple
$(V,D)$ a \textit{V-representation} or vertex representation of
$P$. If the sets $\set{v_i \mid i=1,\dots,k}$ and $\set{d_j \mid
j=1,\dots,\ell}$ are minimal in the sense that no proper subset of any
of the sets generates the same polyhedron, then we call the $v_i$
\textit{vertices} of $P$ and the $d_j$ \textit{extreme directions} of
$P$. The well-known Minkowski-Weyl theorem states that every
polyhedron admits both an H- and a V-representation. Given a
hyperplane $H(w,\gamma)$ with normal vector $w$ and offset $\gamma$, i.e. $H(w,\gamma) =
\set{x \in \R^n \mid w\T x = \gamma}$, we denote the corresponding
halfspace $\set{x \in \R^n \mid w\T x \leq \gamma}$ by
$H^-(w,\gamma)$. Moreover, we say that a hyperplane $H(w,\gamma)$
\textit{supports} a set $C \subseteq \R^n$ or is a \textit{supporting
hyperplane} of $C$, if $C \subseteq H^-(w,\gamma)$ and $C \cap
H(w,\gamma) \neq \emptyset$.

We write $e_i$ for the $i$-th standard unit vector in $\R^n$ and $e$
for the vector whose components are all equal to one.  The \textit{Euclidean
norm} of a vector $x \in \R^n$ is written as $\norm{x}$. The vector
space of real symmetric $n \times n$ matrices is denoted by
$\S^n$. The identity matrix in $\S^n$ is declared $I$. If a matrix $S
\in \S^n$ is positive semidefinite, we write $S \mgeq 0$, and if it is
positive definite, $S \mg 0$. The usual inner product on $\S^n$ is
defined by $S_1 \cdot S_2 = \tr(S_1S_2)$, where $S_1S_2$ means
ordinary matrix multiplication. The \textit{Hausdorff distance}
between two sets $C_1, C_2 \subseteq \R^n$, denoted by
$\haus{C_1}{C_2}$, is defined as
\[
\max \set{\sup_{c_1 \in C_1} \inf_{c_2 \in C_2} \norm{c_1-c_2},
\sup_{c_2 \in C_2} \inf_{c_1 \in C_1} \norm{c_2-c_1} }.
\]
It is well known that the Hausdorff distance defines a metric on the
space of compact subsets of $\R^n$. However, if one argument is
unbounded, the Hausdorff distance may be infinite. In particular, if
$K_1, K_2$ are convex cones, their Hausdorff distance is $0$ if $\cl
K_1 = \cl K_2$ and $+\infty$ otherwise.

The \textit{(linear matrix) pencil} $\A$ \textit{of size $\ell$}
defined by the matrices $A_1,\dots,A_n \in \S^{\ell}$ is the linear
function
\[
\A \colon \R^n \to \S^{\ell}, \quad \A(x) =
\sum_{i=1}^{n} x_iA_i.
\]
A pencil $\A$ of size $\ell$ together with a matrix $A \in
\S^{\ell}$ define the \textit{spectrahedron}
\[
\set{x \in \R^n \mid \A(x) + A \mgeq 0}.
\]
Given two pencils $\A$ and $\B$ of size $\ell$ defined on $\R^n$ and
$\R^m$, respectively, and a matrix $A \in \S^{\ell}$, the set
\[
S := \set{x \in \R^n \mid \exists y \in \R^m \colon \A(x) + \B(y) + A \mgeq
0}
\]
is called a \textit{projected spectrahedron} or \textit{spectrahedral
shadow}. The set $S$ is the projection of the spectrahedron
$\set{(x,y) \in \R^{n+m} \mid (\A+\B)(x,y) + A \mgeq 0}$ onto the
$x$-variables, hence the name.

\section{An approximation algorithm for recession cones of spectrahedra}\label{section3}
In this section we present an algorithm for the computation of a
polyhedral inner and outer approximation of the recession cone of a
spectrahedron. The algorithm is based on augmenting and cutting
schemes introduced in \cite{Kam92} and \cite{BW78}. Augmenting and
cutting schemes are iterative algorithms for compact convex sets that
successively improve the polyhedral approximations. In doing so,
augmenting schemes compute inner approximations and cutting schemes
compute outer approximations. Their functionality can be described as
follows. Assume $C \subseteq \R^n$ is the set to be approximated and
$\O^k$/$\I^k$ are the polyhedral outer and inner approximations
computed after iteration $k$ by a cutting and augmenting scheme,
respectively. Then a cutting scheme refines $\O^k$ by
\begin{enumerate}
\item Choose a direction $w$.
\item Compute a supporting hyperplane $H$ of $C$ with normal vector
$w$.
\item Construct $\O^{k+1} = \O^k \cap H^-$.
\end{enumerate}
In an augmenting scheme the procedure can be stated as
\begin{enumerate}
\item Compute a point $x \in \bd C$.
\item Contruct $\I^{k+1} = \conv \left(\I^k \cup \set{x} \right)$.
\end{enumerate}
These iterative schemes are also jointly called Hausdorff schemes,
because the above steps are typically iterated until the Hausdorff
distance between $C$ and $\O^k$/$\I^k$ is at most a predefined
margin. Various algorithms using the ideas of cutting and augmenting
schemes can be found in the literature, for example to approximate
compact spectrahedral shadows \cite{Cir19} or in the realm of vector
optimization \cite{DLS+21, AUU21}.

As we are working with cones, the Hausdorff distance is unsuitable as
a measure of similarity. Therefore, we consider the following distance
measure on the space of closed convex cones.

\begin{definition}[cf. \cite{IS10}]
Let $K_1, K_2 \subseteq \R^n$ be closed convex cones. The
\textit{truncated Hausdorff distance} between $K_1$ and $K_2$, denoted
$\thaus{K_1}{K_2}$, is defined as
\[
\thaus{K_1}{K_2} := \haus{K_1 \cap B}{K_2 \cap B},
\]
where $B := \set{x \in \R^n \mid \norm{x} \leq 1}$.
\end{definition}

It is known that the truncated Hausdorff distance defines a metric in
the space of closed convex cones in $\R^n$, see \cite{RW98,
Gur62}. Moreover, it clearly holds, that $\thaus{K_1}{K_2} \leq 1$ for
all closed convex cones $K_1$, $K_2$.

We use the truncated Hausdorff distance to define polyhedral
approximations of cones.

\begin{definition}
  Given a closed convex cone $K \subseteq \R^n$ and $\varepsilon > 0$,
  we call a polyhedral cone $R$ an \textit{outer (inner)
  $\varepsilon$-approximation} of $K$, if $K \subseteq R$ ($R
\subseteq K$) and $\thaus{R}{K} \leq \varepsilon$.
\end{definition}

In order to link this idea to the Hausdorff distance and work with
cones in the framework of Hausdorff schemes, we need the concept of
bases.

\begin{definition}\label{def:base}
Let $K \subseteq \R^n$ be a convex cone. A subset $M \subseteq K$ is
called a \textit{base} of $K$, if for every $x \in K\setminus\set{0}$ there exists
a unique $\lambda_x > 0$, such that $\lambda_x x \in M$.
\end{definition}

The following two results establish, that every closed pointed cone
admits a compact base and that such a base can be constructed via its
polar cone. Proofs of these results can be found in \cite[Theorem
3.5]{AT07} and \cite[Theorem 3.2]{Kle57}. There are numerous other
sources that study cones and their properties, e.g., \cite[Chapter
2]{Gre84, BV04, Gru07}.

\begin{proposition}
  Let $K \subseteq \R^n$ be a closed pointed cone. Then $\interior
  K^{\circ} \neq \emptyset$.
\end{proposition}

\begin{proposition}\label{prop:base}
  Let $K \subseteq \R^n$ be a closed pointed cone. Then
  \[
    H(w,-\gamma) \cap K
  \]
  is a compact base of $K$ for ever $w \in \interior K^{\circ}$ and
  $\gamma > 0$.
\end{proposition}

As a consequence, we are able to reduce the approximation of the
recession cone with respect to the truncated Hausdorff distance to the approximation of a compact base of the
cone with respect to the Hausdorff distance.

\begin{proposition}\label{prop:base_apprx}
  Let $K \subseteq \R^n$ be a closed pointed cone and $w \in \interior
  K^{\circ}$. Assume $P$ is a polyhedron satisfying
  $\haus{P}{H(w,-\gamma) \cap K} \leq \varepsilon$. Then the following
  implications hold:
  \begin{enumerate}[label=(\roman*)]
  \item $H(w,-\gamma) \cap K \subseteq P$, $\gamma \geq
    (1+\varepsilon)\norm{w} \Rightarrow \cone P$ is an outer
    $\varepsilon$-approximation of $K$. 
  \item $P \subseteq H(w,-\gamma) \cap K$, $\gamma \geq
    (1+\varepsilon)\norm{w} \Rightarrow \cone P$ is an inner
    $\varepsilon$-approximation of $K$.
  \end{enumerate}
\end{proposition}
\begin{proof}
  We begin the proof with the first assertion. By the Cauchy-Schwarz
  inequality it holds $\norm{x} \geq 1 + \varepsilon$ for every $x
  \in H(w,-\gamma)$. Taking this into account together with
  $\haus{P}{H(w,-\gamma) \cap K} \leq \varepsilon$ yields $\norm{x} \geq
  1$ for every $x \in P$. Proposition \ref{prop:base} ensures that
  $H(w,-\gamma) \cap K$ is a compact set and
    \[
      K = \cone \left(H(w,-\gamma) \cap K\right),
    \]
    because $H(w,-\gamma) \cap K$ is a base of $K$. Since the Hausdorff
    distance between $P$ and $H(w,-\gamma) \cap K$ is bounded from
    above, $P$ must also be compact. This guarantees that $\cone P$ is
    a polyhedral cone.

    Now, assume $\thaus{\cone P}{K}$ is attained as $\norm{p-k}$ for
    $p \in \cone P$ and $k \in K$. This implies, that $\norm{p} =
    1$. Let $\alpha$ be chosen, such that $\alpha k \in P$. By the
    above remark that $\norm{x} \geq 1$ for every $x \in P$ this
    implies $\alpha \geq 1$. Now it holds,
    \begin{equation*}
      \begin{aligned}
        \thaus{\cone P}{K} &= \norm{p-k} \\
        & \leq \alpha \norm{p-k} \\
        & \leq \inf_{\bar{k} \in H(w,-\gamma) \cap K} \norm{\alpha p -
          \bar{k}} \\
        & \leq \haus{P}{H(w,-\gamma) \cap K} \\
        & \leq \varepsilon.
      \end{aligned}
    \end{equation*}
    The second inequality holds true, because $\alpha k$ is the
    projection of $\alpha p$ onto $K$, i.e.
    \[
        \norm{\alpha p - \alpha k} = \inf_{\bar{k} \in K} \norm{\alpha
            p - \bar{k}}.
    \]
    To complete the proof, note that $H(w,-\gamma) \cap K \subseteq P$
    implies $K \subseteq \cone P$, because taking the conical hull is
    an inclusion preserving operation. The proof of the second
    assertion is analogous to the first one with the roles of $P$ and
    $H(w,-\gamma) \cap K$ interchanged.
\end{proof}

Now, from a polyhedral approximation $P$ of a compact base of the
recession cone of a spectrahedron, we obtain an approximation of the
recession cone itself by constructing $\cone P$. Constructing the
conical hull requires no other computation than transposing and
concatenating the matricial data of $P$ as is shown in
\cite[Propositions 8, 10]{CLW19}.

From now on we denote by $C$ the spectrahedron $C := \set{x \in \R^n
  \mid \A(x) + A \mgeq 0}$ for matrices $A, A_1, \dots, A_n \in
\S^{\ell}$. We consider the following semidefinite optimization
problem associated with $C$ and a nonzero direction $w \in \R^n$.
\begin{equation}\tag{P$_1$($\varC,\varw$)}\label{P1}
  \begin{aligned}
    \max \; w\T x \quad \text{s.t.} \quad \A(x)+A \mgeq 0
  \end{aligned}
\end{equation}
Geometrically, solving \eqref{P1} amounts to determining the maximal
shifting of a hyperplane with normal direction $w$ within the
spectrahedron $C$. Another optimization problem we consider reads
as
\begin{equation}\tag{P$_2$($\varC,\varv,\vard$)}\label{P2}
  \begin{aligned}
    \max \; t \quad & \text{s.t.} \quad & \A(x) + A & \mgeq  0 \\
    & & x & = \varv+t\vard,
  \end{aligned}
\end{equation}
for a point $\varv \in \R^n$ and a nonzero direction $\vard \in
\R^n$. Geometrically, given $c \in \interior C$, a solution to
\eqref{P2} is be determined by starting at the point $\varv$ and
moving a positive distance in the direction $\vard$ until a boundary
point of $\varC$ is reached. In general, a part of a solution to
\eqref{P2}, if it exists, will be a point in $\bd C \cap
\set{\varv+t\vard \mid t \in \R}$. In the field of vector
optimization, problems of type \eqref{P2} can be derived from the
\emph{Tammer-Weidner functional}, see \cite{Ger83, GW90}. They also appear
in the literature under the name \emph{Pascoletti-Serafini
scalarization} and are related to Minkowski functionals, see
\cite{PS84, Eic14, Ham06}. The Lagrange dual problem of \eqref{P2} is
\begin{optimizationproblem}{D$_2$($\varC,\varv,\vard$)}{D2}
\min \; (\A(\varv)+A)\cdot U \quad & \text{s.t.} \quad & A_i \cdot U & =  -w_i, \;
i=1,\dots,n \\ & & \vard\T w & = 1 \\ & & U & \mgeq  0.
\end{optimizationproblem}

Solutions to \eqref{P2} and \eqref{D2} give rise to a supporting
hyperplane of $\varC$.

\begin{proposition}\label{prop:p2}
  Let $\varv \in \interior C$ and set $\vard = v-\varv$ for some $v \notin
  C$. Then solutions $(x^*,t^*)$ to \eqref{P2} and $(U^*,w^*)$ to
  \eqref{D2} exist. Moreover, $w^{*\mathsf{T}}x \leq A \cdot U^*$ for
  all $x \in C$ and equality holds for $x=x^*$.
\end{proposition}
\begin{proof}
  Since $\varv \in \interior C$, we can, without loss of generality,
  assume that $\A(\varv) + A \mg 0$, see \cite[Lemma 2.3]{HV07}. Then
  the point $(\varv,0)$ is strictly feasible for \eqref{P2}. Since $v
  \notin C$ and due to convexity of $C$, the first constraint is
  violated whenever $t \geq 1$. The compactness of $C \cap
  \conv\set{c,v}$ implies the existence of a solution $(x^*,t^*)$ of
  \eqref{P2} with $t^* \in [0,1]$. The strict feasibility of the
  problem \eqref{P2} is sufficient for strong duality between the
  problems \eqref{P2} and \eqref{D2} to hold. This implies the
  existence of a solution $(U^*,w^*)$ to \eqref{D2}. Next, let $x \in
  C$ and observe that
  \begin{equation*}
    \begin{aligned}
      A \cdot U^*-w^{*\mathsf{T}}x &= A \cdot U^* + \sum_{i=1}^n
      x_i\left(A_i \cdot U^* \right) \\
      &= A \cdot U^* + \A(x) \cdot U^* \\
      &= \left( \A(x) + A \right) \cdot U^* \\
      &\geq 0.
    \end{aligned}
  \end{equation*}
  Moreover, for $x = x^*$ we have
  \begin{equation*}
    \begin{aligned}
      & w^{*\mathsf{T}}x^* &=& w^{*\mathsf{T}}\left(\varv + t^*d \right)
      \\
      =& - \sum_{i=1}^{n} \left(A_i \cdot U^*\right) c_i
      +t^*w^{*\mathsf{T}}d &=& -\A(\varv) \cdot U^* + t^* \\
      =& -\A(\varv) \cdot U^* + \A(\varv) \cdot U^* + A \cdot U^* &=&
      A \cdot U^*,
    \end{aligned}
  \end{equation*}
  where the fourth equality holds due to strong duality.
\end{proof}

Before we present the algorithm, we make the following assumptions
about $C$:
\begin{enumerate}[label=(C\arabic*)]
\item\label{ass:1} $C$ is line-free.
\item\label{ass:2} There exists $\bar{x} \in \R^n$, such that
  $\A(\bar{x}) \mg 0$.
\end{enumerate}
Since spectrahedra are closed convex sets, Assumption \ref{ass:1}
means that $\recc C$ is pointed.  Moreover, it is straightforward to
verify that $\recc C = \set{x \in \R^n \mid A \mgeq 0}$. Therefore,
Assumption \ref{ass:2} is equivalent to $\interior \recc C \neq
\emptyset$, see \cite[Corollary 5]{RG95}. Pseudo-code of the algorithm
is presentend in Algorithm \ref{alg:1} and one iteration is
illustrated in Figure \ref{fig:alg:1}.

\begin{algorithm}
  \renewcommand\varC{M}
  \DontPrintSemicolon
  \caption{Polyhedral approximation of $\recc C$}\label{alg:1}
  \KwData{Matrix pencil $\A$ and matrix $A$ describing the
    spectrahedron $C$, $\bar{x}$ satisfying $\A(\bar{x}) \mg 0$, error
    tolerance $\varepsilon > 0$}
  \KwResult{An inner $\varepsilon$-approximation $K_{\I}$ and an outer
    $\varepsilon$-approximation $K_{\O}$ of $\recc C$}
  $w \isassigned \left(-A_1 \cdot I, \dots, -A_n \cdot
    I\right)\T/\norm{\left(-A_1 \cdot I, \dots, -A_n \cdot
      I\right)\T}$ \tcp*{interior point of $\left(\recc
      C\right)^{\circ}$} \label{alg1:line1}
  $M \isassigned H^-(w,-(1+\varepsilon)) \cap H^-(-w,1+2\varepsilon)
  \cap \recc C$\; \label{alg1:line2}
  $\bar{x} \isassigned -\frac{2+3\varepsilon}{2w\T \bar{x}}\bar{x}$
  \tcp*{interior point of $M$}\label{alg1:line3}
  $\O, \I \isassigned \emptyset$\; \label{alg1:line4}
  \For{every $z \in \set{-e, e_1, \dots, e_n}$}{\label{alg1:line5}
    Compute a solution $x^*$ to {\renewcommand\varw{z}\eqref{P1}}\; \label{alg1:line6}
    Compute a solution $(y^*,t^*)$ to
    {\renewcommand\varv{\bar{x}}\renewcommand\vard{z}\eqref{P2}}\; \label{alg1:line7}
    $\O \isassigned \O \cap H^-(z,z\T x^*)$\; \label{alg1:line8}
    $\I \isassigned \conv \left(\I \cup \set{y^*}\right)$\; \label{alg1:line9}
  }
  $\kappa \isassigned +\infty$\; \label{alg1:line10}
  \While{$\kappa > \varepsilon$}{
    $\kappa \isassigned \haus{\O}{\I}$\;\label{alg1:line12}
    \For{every vertex $v$ of $\O$}{
      Compute a solution $(U^*, w^*)$ to
      {\renewcommand\varv{\bar{x}}\renewcommand\vard{v-\bar{x}}\eqref{D2}}\; \label{alg1:line14}
      $\O \isassigned \O \cap H^-(w^*, A \cdot U^*)$\; \label{alg1:line15}
    }
    \For{every defining halfspace $H^-(w,\gamma)$ of $\I$}{ 
      Compute a solution $x^*$ to \eqref{P1}\; \label{alg1:line17}
      $\I \isassigned \conv \left(\I \cup \set{x^*}\right)$\; \label{alg1:line18}
    }
  }
  $K_{\O} \isassigned \cone \O$\; \label{alg1:line19}
  $K_{\I} \isassigned \cone \I$\; \label{alg1:line20}
\end{algorithm}

\begin{figure}
  \centering
  \begin{subfigure}[t]{\textwidth}
    \centering
    \includegraphics[]{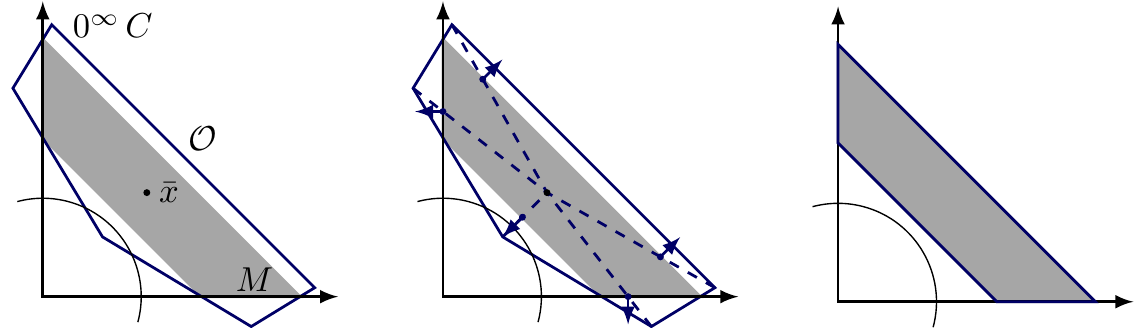}
    \caption{\label{fig:alg:1_outer}\textbf{Left:} The boundary of the
    current outer approximation is colored in blue, the set $M$ is
    shaded in gray. \textbf{Center:} The dotted lines connect the point
    $\bar{x}$ with the vertices of $\O$. Solutions to
    {\renewcommand\varC{M}\renewcommand\varv{\bar{x}}\renewcommand\vard{v-\bar{x}}\eqref{P2}}
    are indicated by the blue dots. The blue arrows are the normal
    vectors of $M$ at these points and correspond to parts of
    solutions to
    {\renewcommand\varC{M}\renewcommand\varv{\bar{x}}\renewcommand\vard{v-\bar{x}}\eqref{D2}}. \textbf{Right:}
    The boundary of the updated outer approximation is shown in blue.}
  \end{subfigure}
  \ \\
  \begin{subfigure}[t]{\textwidth}
    \centering
    \includegraphics[]{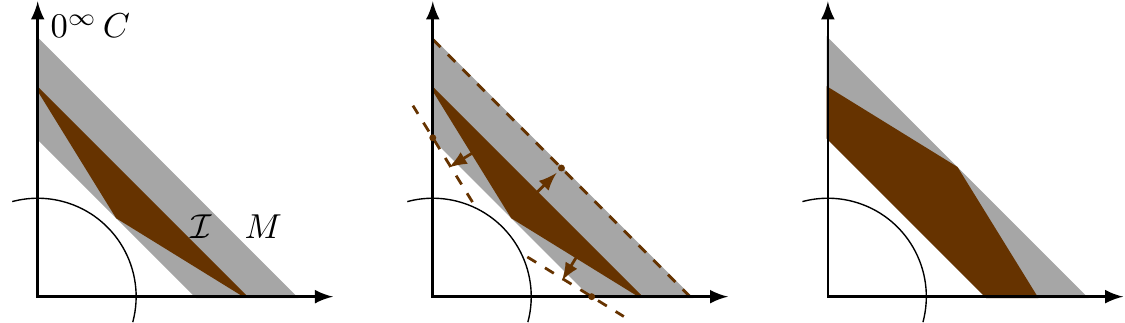}
    \caption{\label{fig:alg:1_inner}\textbf{Left:} The set $M$ is
      shaded in gray, the current inner approximation $\I$ is colored
      orange. \textbf{Center:} The normals of the facets of $\I$ a
      indicated as orange arrows. The dashed lines are translations of
    these facets obtained by solving
    {\renewcommand\varC{M}\eqref{P1}}. The orange dots mark the
    solutions. \textbf{Right:} The updated inner approximation is
    shown in orange.}
  \end{subfigure}
  \decoRule
  \caption{\label{fig:alg:1} One iteration of Algorithm
    \ref{alg:1}. Subfigure \ref{fig:alg:1_outer} illustrates the
    update of the outer approximation and Subfigure
    \ref{fig:alg:1_inner} the update of the inner approximation. The
    boundary of the unit ball is indicated by the black arc to
    illustrate that its intersection with $M$ is empty.}
\end{figure}

The algorithm starts by constructing the set $M$, which can be viewed
as a strip of the cone. $M$ is not a base of $\recc C$ in the sense of
Definition \ref{def:base}, but the set $M \cap H(w,-\gamma)$ is one
for every $\gamma \in [1+\varepsilon,1+2\varepsilon]$. The reason of
working with $M$, instead of with a base, is that $M$ has nonempty
interior, given that $\recc C$ has nonempty interior. This ensures
that Proposition \ref{prop:p2} can be applied. Clearly, from an
approximation of $M$, an approximation of any of the above bases is
available. In line \ref{alg1:line3} the provided interior point
$\bar{x}$ is scaled such that it belongs to the interior of $M$.

In the loop in lines \ref{alg1:line5}--\ref{alg1:line9} an initial outer
and inner approximation of $M$ are computed. The outer approximation is
obtained by solving $n+1$ instances of
{\renewcommand\varC{M}\renewcommand\varw{z}\eqref{P1}}. This gives an
$H$-representation of $\O$. Similarly, the inner approximation is
obtained by solving $n+1$ instances of the problem type
{\renewcommand\varC{M}\renewcommand\varv{\bar{x}}\renewcommand\vard{z}\eqref{P2}}. This
yields $n+1$ points on the boundary of $M$, whose convex hull is a
full dimensional set. Thus, the inner approximation $\I$ is given as a
$V$-representation.

In the main loop the approximations are successively refined according
to the cutting and augmenting procedure introduced at the beginning of
this section. For every vertex $v$ of the outer approximation, one problem
of type {\renewcommand\varC{M}\renewcommand\varv{\bar{x}}\renewcommand\vard{v-\bar{x}}\eqref{D2}}
is solved. Thereby, a supporting hyperplane of $M$ is determined
according to Proposition \ref{prop:p2}. The vertex $v$ is then cut off
from the current outer approximation by intersecting $\O$ with the
corresponding halfspace obtained by
{\renewcommand\varC{M}\renewcommand\varv{\bar{x}}\renewcommand\vard{v-\bar{x}}\eqref{P2}}. To
refine the inner approximation, every defining halfspace or facet of
$\I$ is considered and an instance of
{\renewcommand\varC{M}\eqref{P1}} is solved. This yields a point on
the boundary of $M$, which is appended to the current inner
approximation. The algorithm terminates when the  Hausdorff distance
between the current inner and outer approximation is not larger than
$\varepsilon$. This ensures, that both the Hausdorff distance between
$\O$ and $M$ and the Hausdorff distance between $M$ and $\I$ is at
most $\varepsilon$.

Since $\O$ and $\I$ are polytopes with $\I \subseteq \O$, the quantity
$\haus{\O}{\I}$ will be attained by some vertex of $\O$. Then,
evaluating $\haus{\O}{\I}$ in line \ref{alg1:line12} amounts to
computing $\inf\set{\norm{v-x} \mid x \in \I}$ for every $v \in
\vertices \O$, which are convex quadratic optimization problems, and
taking the maximum over the infima.

Being able to select the vertices of $\O$ requires a
$V$-representation of $\O$. However, during the algorithm $\O$ is
given by an $H$-representation. The task of computing a
$V$-representation from an $H$-representation is called \textit{vertex
  enumeration}. Likewise, in order to choose the facets of $\I$ an
$H$-representation is needed. Computing an $H$-representation from a
$V$-representation is called \textit{facet enumeration}. We point out,
that vertex and facet enumeration are sensitive to the dimension of
the polyhedron and not stable when imprecise arithmetic is used, see
\cite{Ket+08, Loe20}.

\begin{theorem}
  Assume the spectrahedron $C$ input to Algorithm \ref{alg:1}
  satisfies Assumptions \ref{ass:1} and \ref{ass:2}. Then the
  algorithm is finite and works correctly, i.e. it computes an inner
  and outer $\varepsilon$-approximation of $\recc C$.
\end{theorem}
\begin{proof}
  \renewcommand\varC{M}
  Assumption \ref{ass:1} is equivalent to $\recc C$ being
  pointed. This implies that $\polar{\left(\recc C \right)}$ has
  nonempty interior \cite[Corollary 14.6.1]{Roc70} and that the
  matrices $A_1,\dots,A_n$ defining $\A$ are linearly independent, see
  \cite[Lemma 3.2.9]{Net11}. Therefore, the assignment to $w$ in line
  \ref{alg1:line1} is valid. To see, that
  \[
        w = \frac{\left(-A_1 \cdot I, \dots, -A_n \cdot
        I\right)\T}{\norm{\left(-A_1 \cdot I, \dots, -A_n \cdot
        I\right)\T}} \in \interior \polar{\left( \recc{C} \right)},
  \]
  observe that for every $x \in \recc C$, $x \neq 0$, it holds
  \[
    \norm{\left(-A_1 \cdot I, \dots, -A_n \cdot I\right)\T}w\T x =
    -\sum_{i=1}^n x_i\left(A_i \cdot I\right) = - \A(x) \cdot I < 0.
  \]
  The inequality holds, because at least one eigenvalue of $\A(x)$ is
  positive. Now, it follows easily from Proposition \ref{prop:base}
  that the set $M$ defined in line \ref{alg1:line2} is compact. From
  Assumption \ref{ass:2} it follows that the redefinition of $\bar{x}$
  in line \ref{alg1:line3} is an element of $\interior M$. Solutions
  to the problems {\renewcommand\varw{z}\eqref{P1}},
  {\renewcommand\vard{z}\renewcommand\varv{\bar{x}}\eqref{P2}} and
  \eqref{P1} in lines \ref{alg1:line6}, \ref{alg1:line7} and
  \ref{alg1:line17} exist, because $\varC$ is compact. The existence
  of solutions to the problems
  {\renewcommand\varv{\bar{x}}\renewcommand\vard{v-\bar{x}}\eqref{D2}}
  in line \ref{alg1:line14} is guaranteed by Proposition
  \ref{prop:p2}. Note that the component $w^*$ of a solution is always
  nonzero, because $\bar{x} \in \interior M$ implies that the optimal
  value of
  {\renewcommand\varv{\bar{x}}\renewcommand\vard{v-\bar{x}}\eqref{P2}}/{\renewcommand\varv{\bar{x}}\renewcommand\vard{v-\bar{x}}\eqref{D2}}
  is always positive. Now, if the condition of the while loop is
  violated, we have $\haus{\O}{\I} \leq \varepsilon$. Since $\I
  \subseteq \varC \subseteq \O$, this yields $\haus{\O}{\varC} \leq
  \varepsilon$ as well as $\haus{\varC}{\I} \leq \varepsilon$. Taking
  into account the definition of $M$ and the fact $\norm{w} = 1$,
  applying Proposition \ref{prop:base_apprx} yields that $K_{\O}$ is
  an outer $\varepsilon$-approximation of $\recc C$ and $K_{\I}$ is an
  inner $\varepsilon$-approximation of $\recc C$.
  
  The finiteness of the algorithm is due to \cite{Kam93, Kam94}, which
  shows that the inner approximations converge in Hausdorff distance
  to $M$, because they are constructed according to the augmenting
  scheme, and due to \cite[Theorem 4.38]{Cir19}, which shows the
  convergence for the outer approximations.
\end{proof}

\section{An approximation algorithm for recession cones of projections
  of spectrahedra}\label{section4}
In this section, we present an algorithm for the approximation of the
recession cone of a spectrahedral shadow. The technique used in
Algorithm \ref{alg:1} cannot be utilized in this more general
setting. The reason is, that for the approach in Algorithm \ref{alg:1}
it is crucial that an algebraic description of the recession cone is
available, as it gives rise to a description of a base. Unfortunately,
such a description is not easily available in the projected case. As
an example, consider the spectrahedron
\[
  C = \set{\vector[2]{x} \in \R^2 \mid \matrix{x_2 & x_1 \\ x_1 & 1} \mgeq 0}
\]
and its recession cone
\[
  \recc C = \set{\vector[2]{x} \in \R^2 \mid \matrix{x_2 & x_1 \\ x_1 & 0} \mgeq 0}.
\]
$C$ is the epigraph of the function $x \mapsto x^2$ und its recession
cone is simply the cone generated by the direction $(0,1)\T$. Applying
to $C$ the projection $(x_1,x_2)\T \mapsto x_1$ yields the
spectrahedral shadow
\[
  S = \set{x_1 \in \R \mid \exists x_2: \vector[2]{x} \in C}.
\]
Evidently, $S = \R$ and, therefore, $S$ is its own recession
cone. However, $\recc S$ is not equal to the projection of $\recc C$
onto the $x_1$-variable, which results in the singleton $\set{x_1 \in
\R \mid \exists x_2: (x_1,x_2)\T \in \recc C} =
\set{(0,0)\T}$. Consequently, projection is not sufficient to obtain
an algebraic representation of $\recc S$ and Algorithm \ref{alg:1} is
not applicable without one. Therefore, we present in this section an
approximation algorithm that does not directly use information of the
recession cone of a spectrahedral shadow, but only uses information of
the spectrahedral shadow itself to obtain information about its
recession cone. We want to point out however, that it is known, that
the recession cone of a spectrahedral shadow is again a spectrahedral
shadow, see e.g. \cite[Lemma 6.6]{Sch18} and \cite[Proposition
6.1.1]{Net11}.

\begin{remark}
  Clearly, the previous example does not rule out the possibility
  that a semidefinite representation of $\recc S$ can be found. In
  fact, the following relation holds true for closed spectrahedral
  shadows $S$ containing the origin, see \cite[Theorem 14.6]{Roc70}:
  \begin{equation*}
    \recc S = \polar{\left(\cone \polar{S}\right)}.
  \end{equation*}
  Therefore, a representation of $\recc S$ ensues if the conical hull
  and the polar of a spectrahedral shadow can be represented. How one
  can represent the conical hull is certainly well-known, see
  e.g. \cite[Lemma 6.6]{Sch18} and \cite[Proposition 4.3.1]{Net11}. A
  representation of $\polar{S}$ seems to be more difficult to obtain
  or require certain regularity assumptions. In \cite[Proposition
  4.1.7]{Net11} a representation of $\polar{S}$ is presented given
  that the defining pencils of $S$ are strictly feasible. However,
  this condition cannot be guaranteed for the set $\cone
  \polar{S}$. In \cite{Ram97} a representation of the polar of a
  spectrahedron is derived using only the assumption that the set
  contains the origin. This result can be extended to the projected
  case as the polar set of $S$ can be formulated in terms of the polar
  of the spectrahedron from which it is projected, see \cite[Corollary
  16.3.2]{Roc70}. This representation of $\polar{S}$ comes with the
  pitfall that the number of variables needed in the description is
  of the order $O((n+m)\ell^2)$ and the size of the pencil is of the
  order $O((n+m)^2)$, where $n+m$ is the dimension of the
  spectrahedron that is projected and $\ell$ is the size of the
  original pencils. Hence, a representation of $\recc S$ would require
  $O((n+m)^5\ell^2)$ many variables and matrices of size
  $O((n+m)^2\ell^4)$, because the polar of a set has to be computed
  twice according to the equation above.
\end{remark}

Given two pencils $\A$ and $\B$ of equal size $\ell$ and a matrix $A
\in \S^{\ell}$, we consider from now on the spectrahedral shadow
\[
  S = \set{x \in \R^n \mid \exists y \in \R^m : \A(x) + \B(y) + A
    \mgeq 0},
\]
whose recession cone we want to approximate. We make the following
assumptions about $S$:
\begin{enumerate}[label=(S\arabic*)]
\item\label{ass:1proj} $S$ is closed and $S \neq \R^n$.
\item\label{ass:2proj} A point $\bar{x}$ with
  $\A(\bar{x})+\B(\bar{y})+A \mg 0$ for some $\bar{y} \in \R^m$ is
  known. In particular, $\interior S \neq \emptyset$.
\item\label{ass:3proj} A point $\bar{d} \in \interior \recc S \cap B$ is
  known. In particular, $\interior \recc S \neq \emptyset$.
\end{enumerate}

\begin{remark}
  The closedness in Assumption \ref{ass:1proj} cannot be
  omitted. Unlike for the spectrahedral case, where the sets are
  always closed, spectrahedral shadows need not be closed. Consider,
  for example, the set
  \[
    S = \set{x \in \R \mid \exists y \in \R : \matrix{x & 1 \\ 1 & y}
      \mgeq 0},
  \]
  which is equal to the sets $\set{x \in \R \mid \exists y : xy \geq
  1, x \geq 0} = \set{x \in \R \mid x > 0}$. The last set is clearly
  not closed. The closedness is important to guarantee the existence
  of solutions to the problems {\renewcommand\varC{S}\eqref{P1}} and
  {\renewcommand\varC{S}\eqref{P2}} in theory. From a practical
  perspective, however, the fact that $S$ might fail to be closed
  does not affect the ability of modern solvers to compute
  approximate solutions to a prescribed precision, because they
  employ interior-point methods.
\end{remark}

The problems {\renewcommand\varC{S}\eqref{P2} and \eqref{D2}} take the
form
\begin{optimizationproblem}{P$_2$($S,\varv,\vard$)}{P2:proj}
  \max \; t \quad & \text{s.t.} \quad & \A(x) + \B(y) + A & \mgeq 0 \\
  & & x & = \varv + t \vard
\end{optimizationproblem}
and
\begin{optimizationproblem}{D$_2$($S,\varv,\vard$)}{D2:proj}
  \min \; (\A(\varv)+A)\cdot U \quad & \text{s.t.} \quad & A_i \cdot U
  & = -w_i, \; i=1,\dots,n \\ & & B_j \cdot U & = 0, \; j=1,\dots,m \\
  & & \vard\T w & = 1 \\ & & U & \mgeq 0,
\end{optimizationproblem}
where $B_j, j=1,\dots,m$ are the matrices defining the pencil $\B$,
respectively. We derive the following analogous result of Proposition
\ref{prop:p2} for the projected case.

\begin{corollary}\label{cor:p2:proj}
  Let Assumptions \ref{ass:1proj} and \ref{ass:2proj} hold for the
  spectrahedral shadow $S$ and set $d = v-\bar{x}$ for some $v \notin
  S$. Then solutions $(x^*,y^*,t^*)$ to
  {\renewcommand\varv{\bar{x}}\eqref{P2:proj}} and $(U^*,w^*)$ to
  {\renewcommand\varv{\bar{x}}\eqref{D2:proj}} exist. Moreover, $w^*x
  \leq A \cdot U^*$ for all $x \in S$ and equality holds for $x=x^*$.
\end{corollary}
\begin{proof}
  \renewcommand\varv{\bar{x}}
  Assumption \ref{ass:2proj} implies that $\bar{x} \in \interior
  S$. Therefore, the point $(\bar{x},\bar{y},0)$ is strictly feasible for
  \eqref{P2:proj}. Convexity and Assumption \ref{ass:1proj} imply the
  existence of a solution $(x^*,y^*,t^*)$ with $t^* \in [0,1]$. A
  solution $(U^*,w^*)$ of \eqref{D2:proj} exists due to strong
  duality. The rest of the proof is analogous to the proof of
  Proposition \ref{prop:p2} using the fact that $B_j \cdot U = 0$ for
  $j=1,\dots,m$ and feasible points $(U,w)$ of \eqref{D2:proj}.
\end{proof}

We describe the functioning of the algorithm before we present it in
pseudo-code. Similar to Algorithm \ref{alg:1}, the core idea is to
simultaneously maintain an outer approximation $\O$ and an inner
approximation $\I$ of $\recc S$. Also, the approximations are updated
in a similar manner as in Algorithm \ref{alg:1}, i.e. the outer
approximation is updated by computing a supporting hyperplane $H(w,0)$
of $\recc S$ and setting $\O \cap H^-(w,0)$ as the improved
approximation and the inner approximation is updated by computing a
direction $d \in \recc S$ and setting $\cone\left(\I \cup
\set{d}\right)$ as the new inner approximation. Contrary to Algorithm
\ref{alg:1}, however, we cannot steer in advance which approximation
will be updated. After determining a direction $d \in \R^n$, we
investigate
{\renewcommand\varC{S}\renewcommand\varv{\bar{x}}\renewcommand\vard{d}\eqref{P2:proj}}. If
the problem is unbounded, then $d$ is a direction of recession and
$\I$ will be updated. If the problem has an optimal solution, then a
solution $(U^*, w^*)$ to the Lagrange dual yields a supporting
hyperplane $H(w^*,0)$ of $\recc S$ according to Corollary
\ref{cor:p2:proj} and $\O$ is updated. These steps are iterated until
$\thaus{\O}{\I}$ is certifiably not larger than some preset tolerance
$\varepsilon > 0$. To determine initial approximations, we set $\I =
\cone\set{\bar{d}}$ and solve
{\renewcommand\varC{S}\renewcommand\varv{\bar{x}}\renewcommand\vard{-\bar{d}}\eqref{P2:proj}},
from which a supporting hyperplane of $\recc S$ is obtained. The
corresponding halfspace is set as the initial outer approximation. It
remains to explain, how the search directions are chosen in each
iteration. Let $v$ be a vertex of the set $\set{x \in \O \mid
\norm{x}_{\infty} \leq 1}$ and consider the
directions
\[
  d_v^k = \frac{2^k-1}{2^k} v + \frac{1}{2^k} \bar{d}
\]
for $k \in \mathbb{N}$. Starting with $k=1$, we determine whether
$d_v^0 \in \recc S$. If this is the case, then $\I$ is updated and $k$
is incremented. Now, $d_v^1$ is considered. This process is iterated
until either $d_v^k \notin \recc S$ for some $k \in \mathbb{N}$, in
which case $\O$ will be updated and we continue with another vertex
$v$, or $k \geq
\log_2\left(\frac{\norm{v-\bar{d}}}{\varepsilon}\right)$, in which
case $\norm{v-d_v^k} \leq \varepsilon$ and $v$ need no longer be
considerd in the calculations. Note, that
{\renewcommand\varC{S}\renewcommand\varv{\bar{x}}\renewcommand\vard{d_v^k}\eqref{P2:proj}}
only needs to be solved, if $d_v^k$ is not already in $\I$. Otherwise,
it is known, that
{\renewcommand\varC{S}\renewcommand\varv{\bar{x}}\renewcommand\vard{d_v^k}\eqref{P2:proj}}
is unbounded and $k$ can be incremented. The set $\set{x \in \O \mid
\norm{x}_{\infty} \leq 1}$ characterizes $\O$ in the sense that taking
its conical hull yields $\O$ again. In particular, if $v$ is an
extreme direction of $\O$, then $v/\norm{v}_{\infty}$ is a vertex of
$\set{x \in \O \mid \norm{x}_{\infty} \leq 1}$.

\begin{algorithm}
  \renewcommand\varC{S}
  \renewcommand\varv{\bar{x}}
  \renewcommand\vard{\bar{d}}
  \DontPrintSemicolon
  \SetKw{and}{and}
  \caption{Polyhedral approximation of $\recc \varC$}\label{alg:2}
  \KwData{Matrix pencils $\A$, $\B$ and matrix $A$ describing the
    spectrahedral shadow $S$, $\varv$ satisfying \ref{ass:2proj},
    $\bar{d}$ satisfying \ref{ass:3proj}, error tolerance
    $\varepsilon > 0$}
  \KwResult{An inner $\varepsilon$-approximation $\I$ and an outer
    $\varepsilon$-approximation $\O$ of $\recc \varC$}
  Compute a solution $(U^*,w^*)$ to
  {\renewcommand\vard{-\bar{d}}\eqref{D2:proj}}\;
  $\I \isassigned \cone\set{\bar{d}}$\;
  $\O \isassigned H^-(w^*/\norm{w^*},0)$\;
  \Repeat{\upshape{Stop}}{
    Stop $\isassigned$ \texttt{TRUE}\;    
    \For{every nonzero vertex $v$ of $\set{x \in \O \mid \left\lVert x
      \right\rVert_{\infty} \leq 1}$}{\label{line_2:6}
      \For{$k=1$ \KwTo
        $\left\lceil\log_2\left(\frac{\norm{v-\vard}}{\varepsilon}\right)\right\rceil$}{
        $d_v^k \isassigned \frac{2^k-1}{2^k} v + \frac{1}{2^k}
        \bar{d}$\;
        \If{$d_v^k \in \I$}{
          continue\;
        }\ElseIf{\renewcommand\vard{d_v^k}\eqref{P2:proj} is
          {\normalfont\texttt{UNBOUNDED}}}{
          $\I \isassigned \cone\left(\I \cup \set{d_v^k}\right)$\;
        }\Else{
          Compute a solution $(U^*,w^*)$ to
          {\renewcommand\vard{d_v^k}\eqref{D2:proj}}\;\label{line_2:14}
          $\O \isassigned \O \cap H^-(w^*/\norm{w^*},0)$\;
          Stop $\isassigned$ \texttt{FALSE}\;
          break\;\label{line_2:17}
        }
      }
    }
  }
\end{algorithm}

\begin{figure}
  \renewcommand\varC{S}
  \renewcommand\varv{\bar{x}}
  \centering
  \includegraphics[]{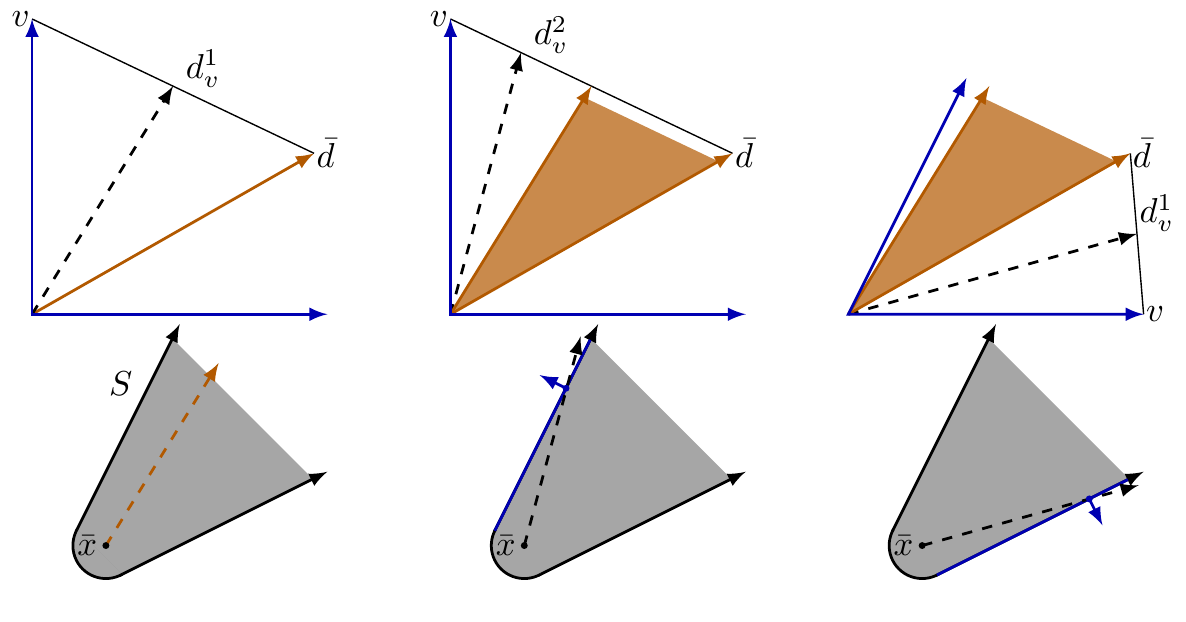}
  \decoRule
  \caption{\label{fig:alg:2} One iteration of the outer for loop of
    Algorithm \ref{alg:2} for the direction $v$. The top row shows how
    the blue outer approximation and the orange inner approximation
    change during the process and how the next search direction is
    obtained. The bottom row depicts the set $S$ with interior point
    $\bar{x}$ and illustrates the problems
    {\protect\renewcommand\vard{d_v^k}\eqref{P2:proj}/\eqref{D2:proj}} that are
    considered. In the leftmost column the direction $d_v^1$ is chosen
    as the midpoint of the line segment from $v$ to $\bar{d}$ and
    {\protect\renewcommand\vard{d_v^1}\eqref{P2:proj}} is determined to be
    unbounded as indicated by the dashed orange ray originating from
    $\bar{x}$. Thus, the inner approximation is updated by
    constructing the conical hull of $\bar{d}$ and $d_v^1$. This is
    shown in the top of the center column. Since $d_v^1$ is a
    direction of recession, $d_v^2$ is obtained as the midpoint
    between $v$ and $d_v^1$. The bottom picture illustrates part of
    the solution of {\protect\renewcommand\vard{d_v^2}\eqref{D2:proj}} as the
    outward pointing normal direction, which yields a supporting
    hyperplane of $S$. Subsequently, the outer approximation is
    updated by shifting the hyperplane to the origin and intersecting
    the outer approximation with the corresponding halfspace. The
    right column shows the updated outer approximation. Since $\O$ was
    updated, a new direction $v$ is chosen.}
\end{figure}

One iteration of Algorithm \ref{alg:2} for one direction $v$ is
illustrated in Figure \ref{fig:alg:2}. The outer for loop in line
\ref{line_2:6} requires a V-representation of the polytope $\set{x \in
\O \mid \norm{x}_{\infty} \leq 1}$. During the execution of the
algorithm, however, $\O$ is always given by an H-representation and
$\I$ by a V-representation. Therefore, a vertex enumeration has to be
performed in every iteration. The algorithm terminates if lines
\ref{line_2:14} - \ref{line_2:17} are not executed during one full
iteration, i.e. if $\O$ receives no upgrade. This happens only if for
every $v$ of the set in line \ref{line_2:6} and $k =
\left\lceil\log_2\left(\frac{\norm{v-\vard}}{\varepsilon}\right)\right\rceil$
the point $d_v^k$ belongs to $\I$. This implies
\[
\inf_{d_i \in \I} \norm{v-d_i} \leq \norm{v-d_v^k} =
\frac{1}{2^k}\norm{v-\bar{d}} \leq \varepsilon
\]
for every vertex $v$, which certifies $\thaus{\O}{\I} \leq
\varepsilon$. Since $\I \subseteq \recc S \subseteq \O$, the relations
$\thaus{\O}{\recc S} \leq \varepsilon$ and $\thaus{\recc S}{\I} \leq
\varepsilon$ follow. We now prove that Algorithm \ref{alg:2} works
correctly and is finite.

\begin{theorem}
  Let Assumptions \ref{ass:1proj}-\ref{ass:3proj} be satisfied for a
  spectrahedral shadow $S$. Then Algorithm \ref{alg:2} terminates after
  finitely many steps and is correct, i.e. the algorithm returns an
  inner $\varepsilon$-approximation $\I$ and an outer
  $\varepsilon$-approximation $\O$ of $\recc S$.
\end{theorem}
\begin{proof}
  \renewcommand\varC{S}\renewcommand\varv{\bar{x}}
  We first prove the correctness. Assumptions \ref{ass:1proj} and
  \ref{ass:3proj} imply that $-\bar{d} \notin \recc S$. Otherwise we
  had $0 \in \interior \recc S$, which implies $S = \R^n$, see
  \cite[Theorem 6.1]{Roc70}. Then, since $\A(\bar{x})+\B(\bar{y})+A
  \mg 0$ and $S$ is closed, a solution to
  {\renewcommand\vard{-\bar{d}}\eqref{D2:proj}} in line 1 exists
  according to Corollary \ref{cor:p2:proj}. Furthermore, Assumption
  \ref{ass:3proj} and Corollary \ref{cor:p2:proj} imply that the set
  $\I$ initialized in line 2 and the set $\O$ initialized in line 3
  are an inner and an outer approximation of $\recc S$,
  respectively. For the latter, take into account that the recession
  cone of a polyhedron in H-representation is described by the
  corresponding homogeneous system of inequalities, see
  \cite[p. 62]{Roc70}. During the first execution of the outer for
  loop in lines 6--17, $\O = H^-(w^*/\norm{w^*},0)$ and the set
  $\set{x \in \O \mid \left\lVert x \right\rVert_{\infty} \leq 1}$ in
  line \ref{line_2:6} is nonempty. In particular, it is a polytope and the for loop
  is executed at least once. Now, let $v$ be a from this set, for
  which the inner for loop is executed. Without loss of generality we
  can assume that such $v$ exists, otherwise we have
  \[
    \log_2\left(\frac{\norm{v-\bar{d}}}{\varepsilon}\right) \leq 0
  \]
  for all vertices $v$ of $\set{x \in \O \mid \left\lVert x
  \right\rVert_{\infty} \leq 1}$. This implies that
  \begin{align*}
    \thaus{\O}{\I} &= \sup_{v \in \O \cap B} \inf_{d \in \I}
    \norm{v-d}\\
    &\leq \sup_{v \in \O \cap \set{x \mid \norm{x}_{\infty} \leq 1}}
    \inf_{d \in \I} \norm{v-d}\\
    &\leq \max\set{\norm{v-\bar{d}} \mid v \text{ is a vertex of } \O
    \cap \set{x \mid \norm{x}_{\infty} \leq 1}}\\
    & \leq \varepsilon.
  \end{align*}
  The first equality holds true because $\I \subseteq \O$ and the
  projection onto a convex cone is a norm-reducing operation, i.e. the
  infimum is attained for some $d \in B$, cf. \cite{Mor62}. The
  inequality in the second line simply reflects the fact that $B$ is
  contained in $\set{x \in \R^n \mid \norm{x}_{\infty} \leq
  1}$. Moreover, the supremum in the second line is attained at a
  vertex of $\O \cap \set{x \in \R^n \mid \norm{x}_{\infty} \leq 1}$,
  see \cite[Theorem 3.3]{Bat86}. Thus, $\O$ and $\I$ are
  $\varepsilon$-approximations of $\recc S$ already.

  Next, note that the search direction $d_v^k$ defined in line 8 is
  zero only if $v = -\bar{d}$. However, due to the fact that $\bar{x} \in
  \interior \S$ it is straightforward to see, that the initial outer
  approximation in line 3 does not contain $-\bar{d}$. Hence, subsequent
  approximations also do not contain $-\bar{d}$, i.e. $d_v^k \neq 0$ for
  every vertex $v$ selected in line \ref{line_2:6}. Now, assume
  {\renewcommand\vard{d_v^k}\eqref{P2:proj}} is unbounded for some $k$,
  i.e. $\bar{x} + td_v^k \in S$ for every $t \geq 0$. Since $S$ is
  closed by Assumption \ref{ass:1proj}, this implies $d_v^k \in \recc
  S$, see \cite[Theorem 8.3]{Roc70}. Hence, the updated inner
  approximation in line 12 still satisfies $\I \subseteq \recc S$. If
  {\renewcommand\vard{d_v^k}\eqref{P2:proj}} is not unbounded, then
  Assumption \ref{ass:2proj} and Corollary \ref{cor:p2:proj} guarantee
  that {\renewcommand\vard{d_v^k}\eqref{D2:proj}} has an optimal solution
  $(U^*,w^*)$. Moreover, $H(w^*,A \cdot U^*)$ supports $S$, which
  implies that $H(w^*/\norm{w^*},0)$ supports $\recc S$. Therefore,
  after the update in line 15 it still holds $\recc S \subseteq
  \O$ and $\O$ is a cone.

  It remains to show that upon termination of the algorihm, $\I$ and
  $\O$ are $\varepsilon$-approximations of $\recc S$. Algorithm 2
  terminates if and only if during one iteration of the repeat until
  loop in lines 4--18 the outer approximation is not modified,
  i.e. lines 14--17 are not executed. This is the case if and only if
  for vertex $v$ of $\set{x \in \O \mid \norm{x}_{\infty} \leq 1}$ it
  holds $d_v^{\bar{k}(v)} \in \I$ with $\bar{k}(v) =
  \left\lceil\log_2\left(\frac{\norm{v-\vard}}{\varepsilon}\right)\right\rceil$.
  Using the same derivation as above but replaing $\bar{d}$ with
  $d_v^{\bar{k}(v)}$ and taking into accound the relation
  \begin{equation*}
    \begin{aligned}
      \norm{v-d_v^{\bar{k}(v)}} &\leq
      \norm{v-\frac{\norm{v-\vard}-\varepsilon}{\norm{v-\vard}}v-\frac{\varepsilon}{\norm{v-\vard}}\vard}
      \\
      &= \frac{\varepsilon}{\norm{v-\vard}}\norm{v-\vard} \\
      &= \varepsilon
    \end{aligned}
  \end{equation*}
  this yields
  \[
    \thaus{\O}{\I} \leq \varepsilon.
  \]
  Since $\I \subseteq \recc S \subseteq \O$ this implies both
  $\thaus{\O}{\recc S} \leq \varepsilon$ and $\thaus{\recc S}{\I} \leq
  \varepsilon$. This concludes the proof of correctness.    
  
  Now, we show that Algorithm \ref{alg:2} indeed terminates after
  finitely many iterations. Assume that at some point during the
  algorithm the approximations $\O$ and $\I$ have been computed and
  that the algorithm does not halt for these sets. Then there exists
  a vertex $v$ of the set in line \ref{line_2:6} and some $k \leq
  \left\lceil\log_2\left(\frac{\norm{v-\bar{d}}}{\varepsilon}\right)\right\rceil$
  for which {\renewcommand\vard{d_v^k}\eqref{P2:proj}} is bounded. This
  implies that lines 14--17 are executed. Let ($U^*,w^*$) be a
  solution to {\renewcommand\vard{d_v^k}\eqref{D2:proj}} and let $w$ be the
  normalized vector $w^*/\norm{w^*}$. Then, by the same reasoning as
  above, the hyplerplane $H(w,0)$ supports $\recc S$. Therefore, $\I
  \subseteq H^-(w,0)$. By Assumption \ref{ass:3proj} there exists some
  $r>0$, such that $B_r(\bar{d}) := \set{x \in \R^n \mid
    \norm{x-\bar{d}} \leq r} \subseteq \I$. Hence,
  \begin{equation*}
    \begin{aligned}
      \inf\set{\norm{x-\vard} \mid x \in H(w,0)} &= \left\lvert
        w\T \vard \right\rvert \\
      &= -w\T\vard \\
      &\geq r.
    \end{aligned}  
  \end{equation*}
  Taking into account the definition of $d_v^k$ and the fact that
  $w^{*\mathsf{T}}d_v^k = 1$, we conclude
  \begin{equation*}
    \begin{aligned}
      \inf\set{\norm{x-v} \mid x \in H(w,0)} &= \left\lvert w\T v
      \right\rvert \\
      &= w\T v \\
      &\geq w\T
      \left(\frac{2^k}{2^k-1}d_v^k-\frac{1}{2^k-1}\bar{d} \right) \\
      &= \frac{2^k}{(2^k-1)\norm{w^*}}+\frac{r}{2^k-1} \\
      &\geq \frac{r}{2^k-1} \\
      &\geq \frac{r}{\frac{2\norm{v-\bar{d}}}{\varepsilon}-1} \\
      &> \frac{\varepsilon r}{4\sqrt{n}-\varepsilon}.
    \end{aligned}
  \end{equation*}
  The last inequality follows from $\norm{v-\bar{d}} <
  2\sqrt{n}$. That means, whenever $\O$ receives an update in line 15,
  a ball of radius at least $\frac{\varepsilon
  r}{4\sqrt{n}-\varepsilon}$ around $v$ is cut off. Since $v$ belongs
  to the compact set $\set{x \in \R^n \mid \norm{x}_{\infty} \leq 1}$,
  $\O$ can only be updated a finite number of times. This proves that
  the algorithm terminates.
\end{proof}

In comparison to Algorithm \ref{alg:1}, the main advantage of
Algorithm \ref{alg:2} is its higher degree of flexibility. It
generalizes Algorithm \ref{alg:1} in two aspects. Firstly it works for
spectrahedral shadows, which contain spectrahedra as a proper
subclass, and secondly it is also applicable to sets with a recession
cone that is not necessarily pointed. The only drawback of Algorithm
\ref{alg:2} is that an interior point of the recession cone must be
known beforehand and it is not trivial to compute one. However,
knowledge of an interior point of the set to be approximated is a
commom assumption in the literature about polyhedral approximation,
see e.g. \cite{Kam92}, where it is assumed, that the origin is
contained in the interior of the set, \cite{DG86}, where some form of
Slater's constraint qualification is assumed, or \cite{AUU21}, where
some interior point is assumed to exist.

\section{Examples}\label{section5}
We apply Algorithms \ref{alg:1} and \ref{alg:2} to three examples in
this section and provide numerical results. The algorithms are
implemented in Python version 3.9.12. The semidefinite problems
(P$_1$), (P$_2$) and (D$_2$) are solved using version 1.1.18 of the
API \texttt{CVXPY} \cite{DB16,AVD+18} with the solver \texttt{SCS}
\cite{OCP+21}. To carry out the vertex and facet enumerations on the
approximations we use the Python wrapper \texttt{pycddlib} of the
library \texttt{cddlib} \cite{FP96}. The figures in Examples
\ref{ex:2} and \ref{ex:3} are generated with the polyhedral
calculus toolbox \texttt{bensolve tools} \cite{LW16,CLW19}. The first
example considers a spectrahedron and is applied to both algorithms,
while in Examples \ref{ex:2} and \ref{ex:3} spectrahedral shadows are
investigated where only Algorithm \ref{alg:2} can be applied.
\begin{example}\label{ex:1}
  We consider the spectrahedron
  \[
    C = \set{x \in \R^3 \mid \matrix{x_1 & x_2 \\ x_2 & x_3} \mgeq
      \matrix{1 & 0 \\ 0 & 1}}.
  \]
  Its recession cone is the set $S_+^2$. As input parameters we set
  $\bar{x} = \frac{\sqrt{2}}{2}(1,0,1)\T$ for Algorithm \ref{alg:1} and
  $\bar{x} = (2,0,2)\T$, $\bar{d} = \frac{\sqrt{2}}{2}(1,0,1)\T$ for
  Algorithm \ref{alg:2}. Figure \ref{fig:ex:1} shows the number of
  solved SDP subproblems and the elapsed CPU time for different values
  of $\varepsilon$ and both algorithms. We see that the number of solved
  subproblems is of the same order of magnitude for both
  algorithms. Hoewever, Algorithm \ref{alg:2} performs the approximation
  much faster, by a factor of around 23 for $\varepsilon=10^-3$. This is
  explained by the fact that the number of solved SDPs for Algorithm
  \ref{alg:2} includes the unbounded instances of
  {\renewcommand\varv{\bar{x}}\renewcommand\vard{\bar{d}}\eqref{P2}},
  for which \texttt{SCS} can terminate as soon as the unboundedness is
  detected. Another remarkable observation is that the number of
  subproblems for Algorithm \ref{alg:1} is smaller for
  $\varepsilon=0.05$ than for $\varepsilon=0.1$. The reason is that
  the set $M$ approximated by Algorithm \ref{alg:1} depends on the
  value of $\varepsilon$. Recall that $M$ can be seen as a slice of
  $\recc C$ and as $\varepsilon$ decreases so does the thickness of
  $M$. This evidently has an impact on the performance.

  \begin{figure}
    \centering
    \includegraphics[]{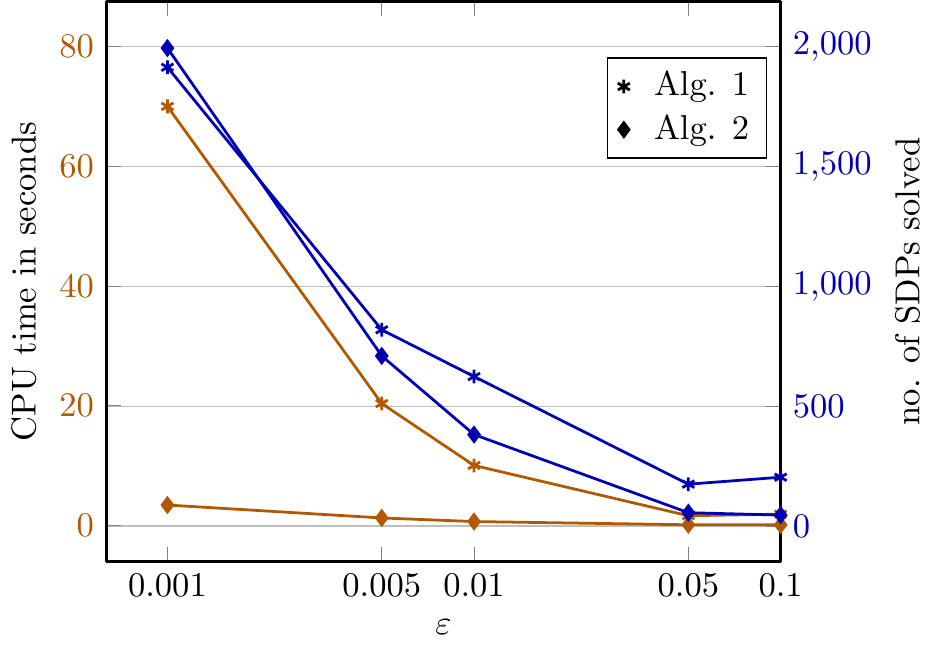}
    \decoRule
    \caption{\label{fig:ex:1}Computational results for Example
      \ref{ex:1}.}
  \end{figure}
  
\end{example}

\begin{example}\label{ex:2}
  We consider the cone $\Sigma_{1,4}$ of univariate polynomials of
  degree at most four which can be represented as a sum of squares of
  polynomials. By identifying the polynomials by their coefficients it
  can be written as the following spectrahedral shadow in $\R^5$:
  \[
    \Sigma_{1,4} = \set{x \in \R^5 \mid \exists y \in \R: \matrix{x_1 & \frac{1}{2} x_2
        & \frac{1}{3} x_3-y \\ \frac{1}{2} x_2 & \frac{1}{3} x_3+2y &
        \frac{1}{2} x_4 \\ \frac{1}{3} x_3-y & \frac{1}{2} x_4 & x_5}
      \mgeq 0}.
  \]
  The cone of sums of squares of polynomials or SOS cone is contained
  in the cone of nonnegative polynomials. It is of interest, because
  optimization problems involving SOS polynomials are tractable by
  semidefinite programming \cite{PT20}. For a derivation of
  semidefinite representations of SOS polynomials see
  \cite[p.~61]{BPT13}. We approximate $\Sigma_{1,4}$ with $\varepsilon=0.1$
  and $\bar{x}=\bar{d}=1/\sqrt{5}e$. Algorithm \ref{alg:2} computes
  the approximations in $26.4$ seconds while solving $1081$
  subproblems. Figure \ref{fig:ex:2} shows the projections of the base
  \[
    \Sigma_{1,4} \cap \set{x \in \R^5 \mid p\T x=-1},
  \]
  where $p \approx (-0.71,-0.08,-0.21,-0.08,-0.66)\T$,
  onto the variables
  \[
    \matrix{x_1 \\ x_2 \\ x_3}, \matrix{x_1 \\ x_2 \\ x_4},
    \matrix{x_1 \\ x_3 \\ x_4} \text{ and } \matrix{x_2 \\ x_3 \\ x_4},
  \]
  respectively.
  
  \begin{figure}
    \centering
    \begin{subfigure}[t]{0.45\textwidth}
      \centering
      \includegraphics[scale=0.2]{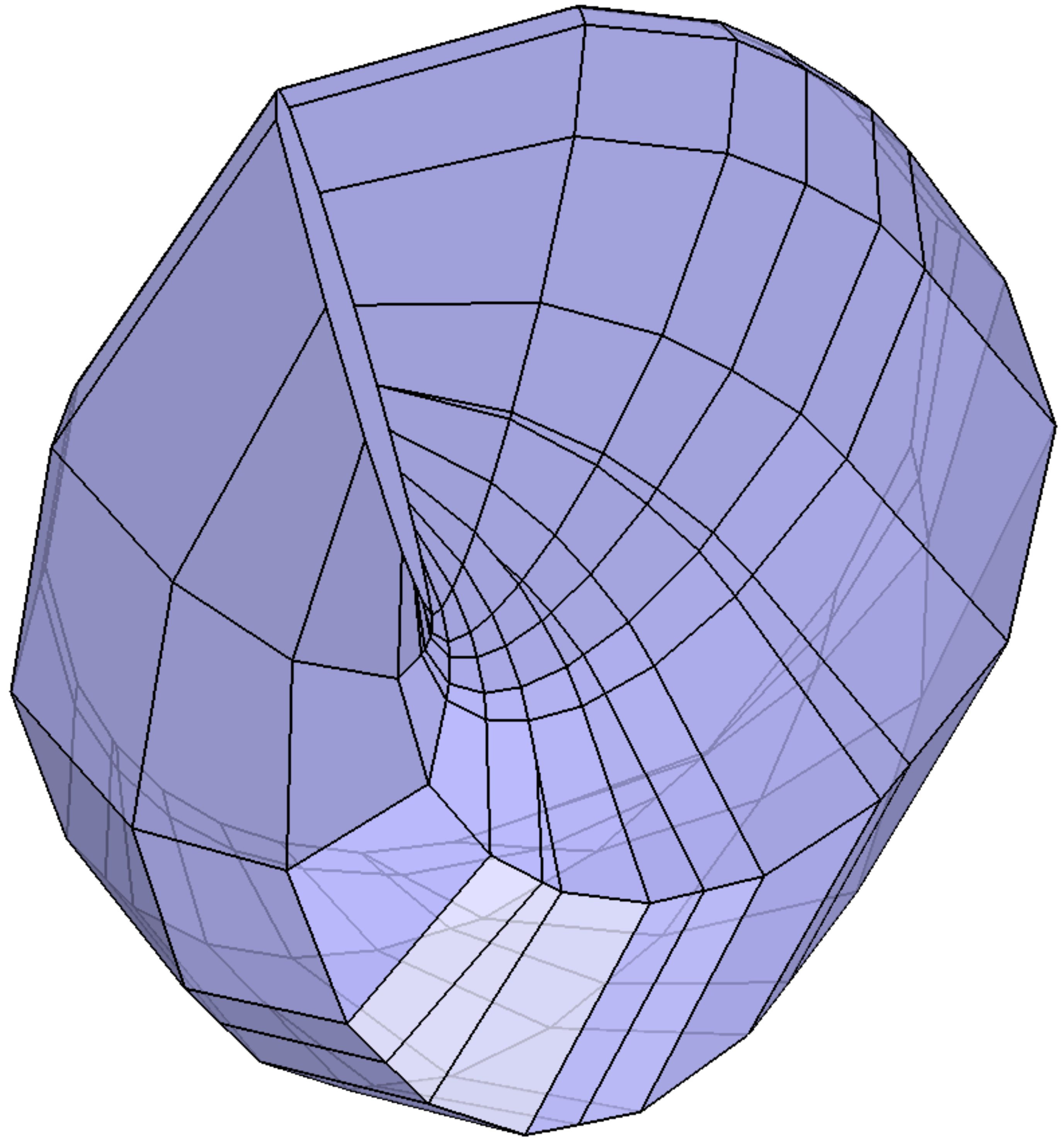}
    \end{subfigure}
    \hspace{.5cm}
    \begin{subfigure}[t]{0.45\textwidth}
      \centering
      \includegraphics[scale=0.2]{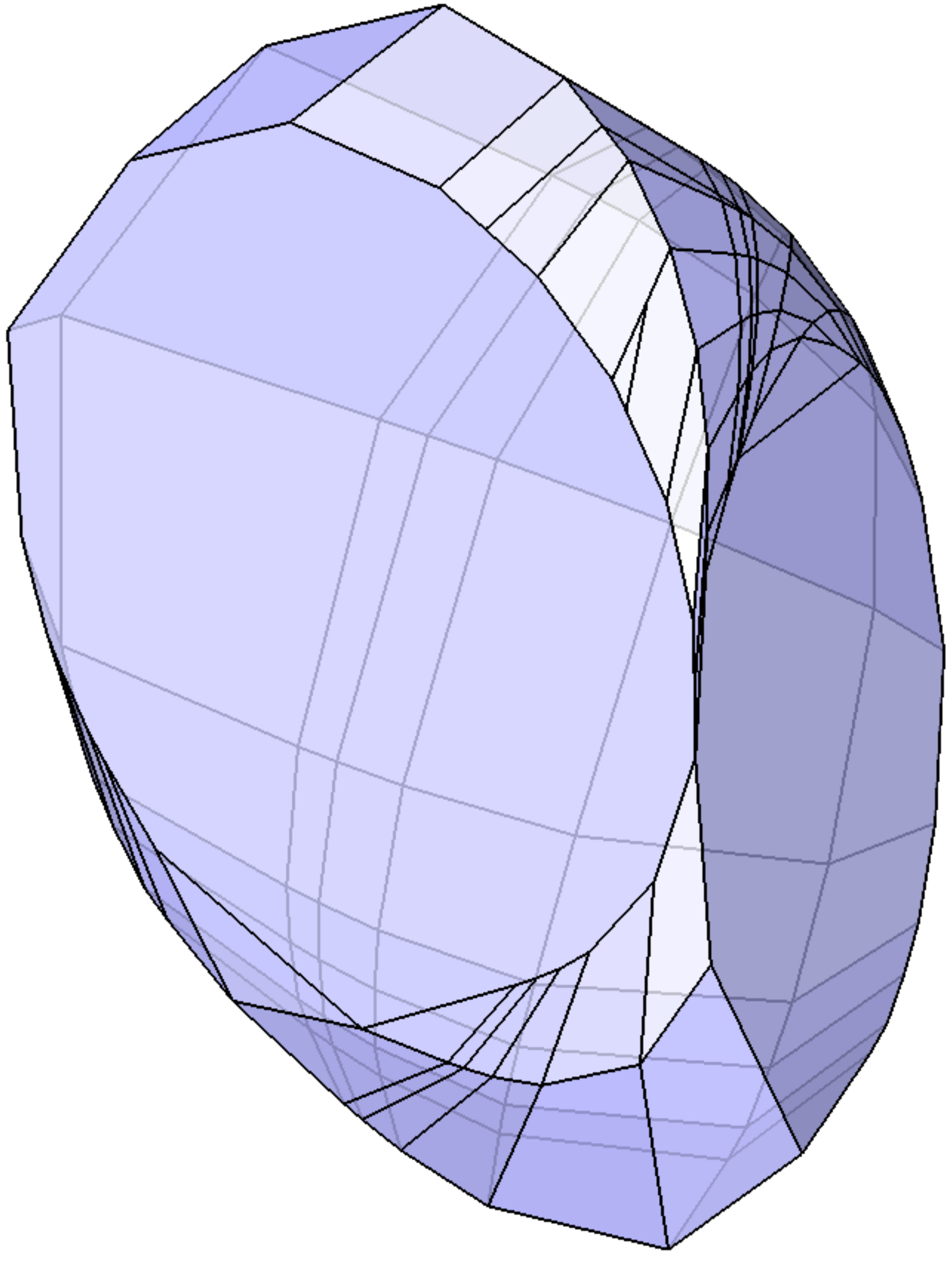}
    \end{subfigure}
    \ \\
    \ \\
    \begin{subfigure}[t]{0.45\textwidth}
      \centering
      \includegraphics[scale=0.17]{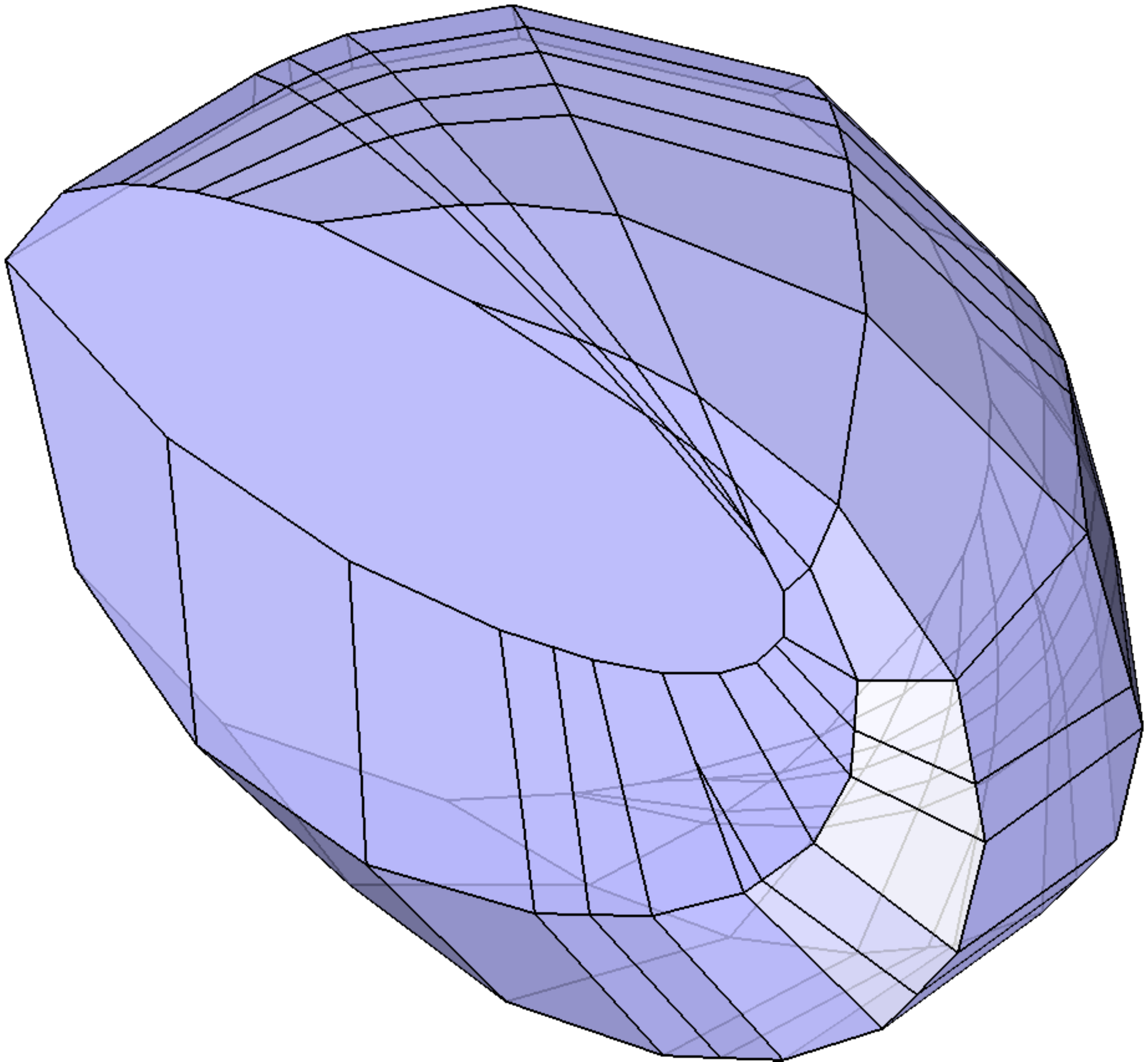}
    \end{subfigure}
    \hspace{.5cm}
    \begin{subfigure}[t]{0.45\textwidth}
      \centering
      \includegraphics[scale=0.2]{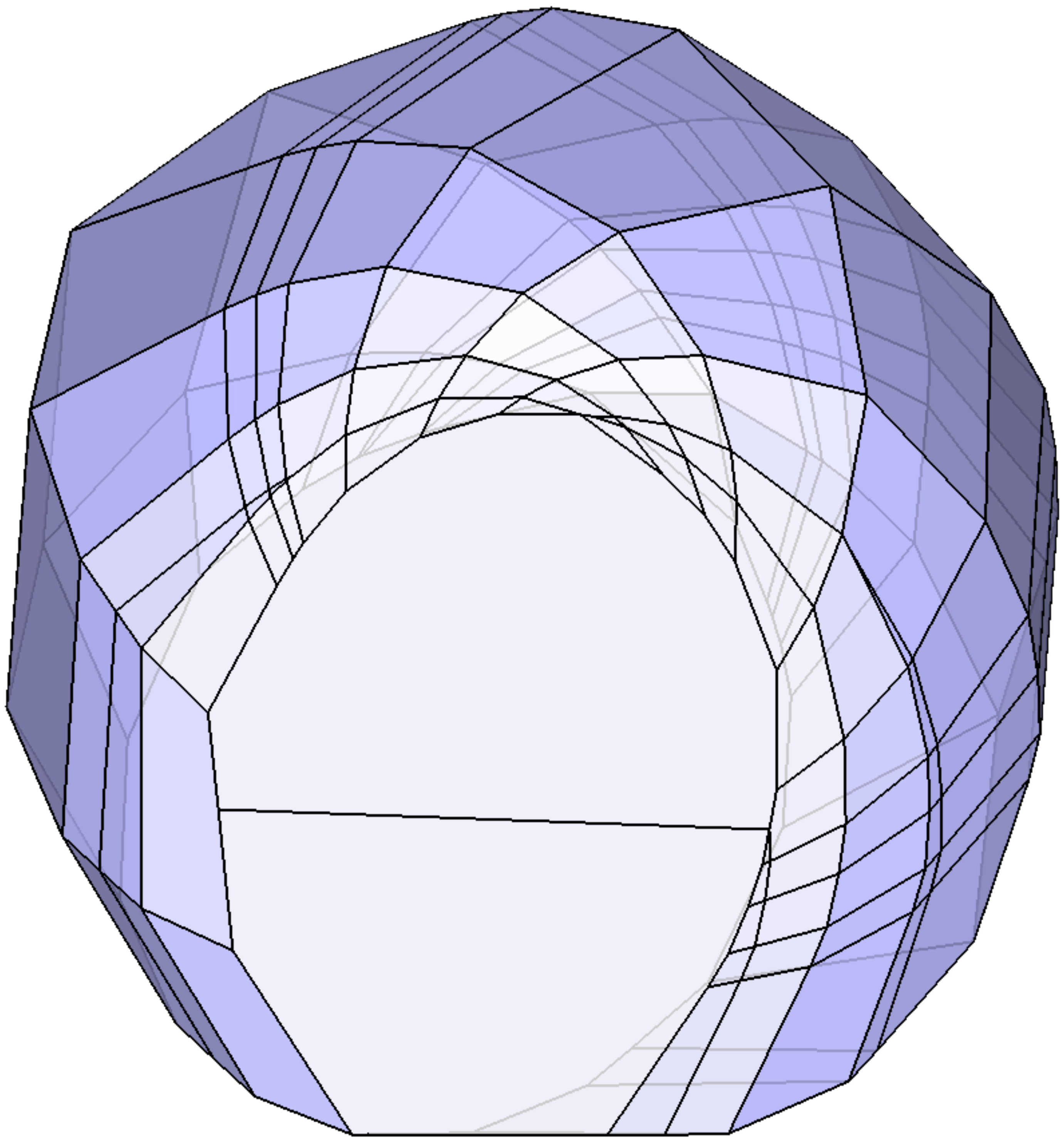}
    \end{subfigure}
    \decoRule
    \caption{\label{fig:ex:2} The orthogonal projections of the
      4-dimensional base of the cone $\Sigma_{1,4}$ from Example
      \ref{ex:2} onto three coordinates.}
  \end{figure}
  
\end{example}

\begin{example}\label{ex:3}
  Since Algorithm \ref{alg:2} operates in the ambient space of the
  spectrahedral shadow and not in the space of the typically high
  dimensional spectrahedron, we want to investigate how well it scales
  with the size of the involved pencils. Therefore, we consider the
  shadow
  \[
    S = \set{\matrix{x_{11} \\ x_{22} \\ x_{33} \\ \ell(X)} \mid X
      \mgeq I},
  \]
  where $\ell(X) = x_{11} + x_{22} + x_{33} + 2(x_{12} + x_{13} +
  x_{23})$, for different sizes $n$ of the matrix $X$. Since $X$ is
  symmetric the number of variables depends on $n$ as
  $\frac{n(n+1)}{2}$. For the computations we set $\bar{d} =
  (1,1,1,0)\T$. Figure \ref{fig:ex:3} shows the outer and inner
  approximations of the base
  \[
    H\left(-\frac{1}{\sqrt{4}}e,-1\right) \cap \recc{S}
  \]
  for $n=3$ and $\varepsilon = 0.01$, respectively. The set is the
  convex dual of a spectrahedron known as \emph{elliptope}
  \cite{MS21}, which has applications in statistics. Some numerical
  results are presented in Table \ref{tbl:ex:3}. It shows the size $n$
  of the pencil, the numver of variables (\#var), the number of solvex
  SDPs (\#SDP) and the elapsed time. The number of solved SDPs for
  $\varepsilon=0.01$ is approximately the same for different $n$. For
  $n \geq 7$ the vertex enumeration fails and no results are
  obtainable. For $\varepsilon=0.1$ the computation time generally
  increases, although drops can be seen, for example for $n=9$. This
  may be explained by \texttt{SCS} taking advantage of the problem
  structure that is hidden from the user.
  
  \begin{figure}
    \centering
    \begin{subfigure}[t]{.45\textwidth}
      \centering
      \includegraphics[scale=0.51]{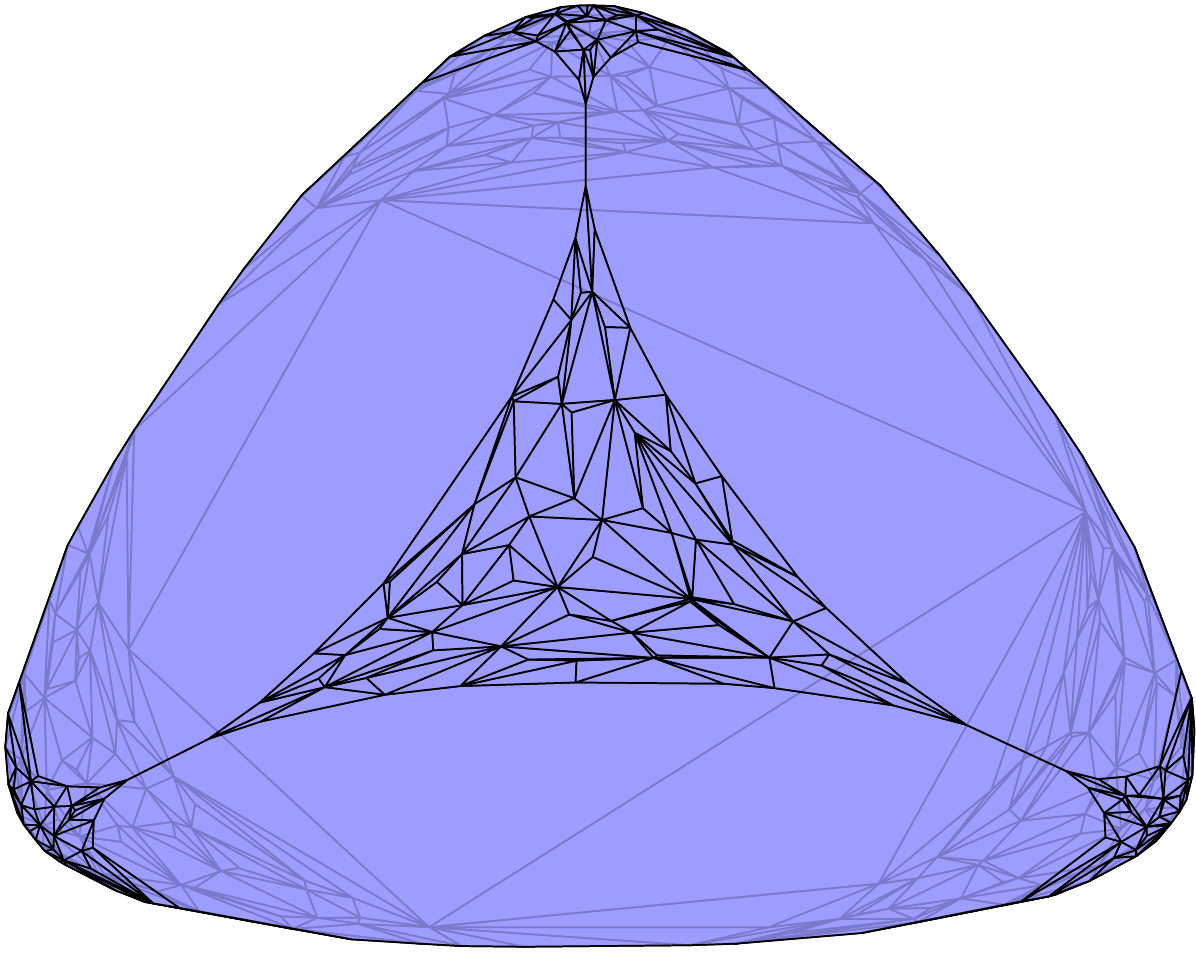}
    \end{subfigure}
    \hspace{1cm}
    \begin{subfigure}[t]{.45\textwidth}
      \includegraphics[scale=0.4]{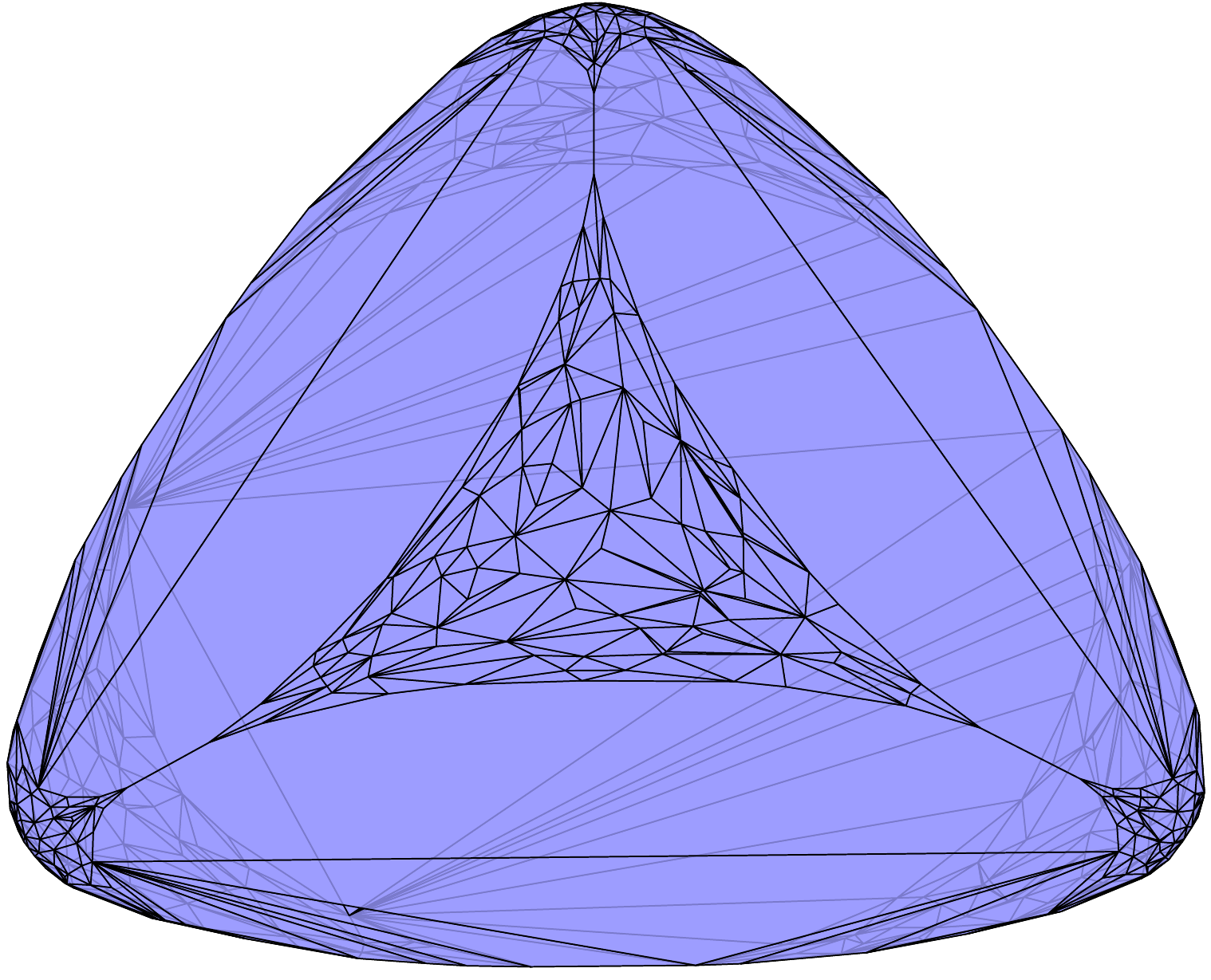}
    \end{subfigure}
    \decoRule
    \caption{\label{fig:ex:3} The outer and inner approximations of a
      base of $\recc{S}$ computed by Algorithm \ref{alg:2},
      respectively.}
  \end{figure}

  \begin{table}
    \centering
    \begin{subtable}[t]{.45\textwidth}
      \centering
      \begin{tabular}{rrrr}
        \toprule
        \boldmath{$n$} & \textbf{\#var} & \textbf{\#SDP} & \textbf{time}
        \\
        \midrule
        3 & 6 & 474 & 2.25 \\
        6 & 21 & 673 & 12.92 \\
        9 & 45 & 384 & 6.10 \\
        12 & 78 & 673 & 23.84 \\
        15 & 120 & 721 & 45.96 \\
        \bottomrule
      \end{tabular}
    \end{subtable}
    \hspace{.5cm}
    \begin{subtable}[t]{.45\textwidth}
      \centering
      \begin{tabular}{rrrr}
        \toprule
        \boldmath{$n$} & \textbf{\#var} & \textbf{\#SDP} & \textbf{time}
        \\
        \midrule
        3 & 6 & 11810 & 132.65 \\
        4 & 10 & 11879 & 111.39 \\
        5 & 15 & 11647 & 152.99 \\
        6 & 21 & 12059 & 308.07 \\
        \bottomrule
      \end{tabular}
    \end{subtable}
    \caption{\label{tbl:ex:3} Numerical results for Example \ref{ex:3}
    for $\varepsilon=0.1$ and $\varepsilon=0.01$, respectively.}
  \end{table}
  
\end{example}
 
\bibliographystyle{plain}
\bibliography{bib/references}

\end{document}